\documentclass{amsart}
\usepackage{amsmath}
\usepackage{amssymb}
\newtheorem{thm}{Theorem}[section]
\newtheorem{cor}[thm]{Corollary}
\newtheorem{conj}[thm]{Conjecture}
\newtheorem{lem}[thm]{Lemma}
\newtheorem{prop}[thm]{Proposition}
\newtheorem{cons}[thm]{Construction}
\theoremstyle{definition}
\newtheorem{defn}[thm]{Definition}
\theoremstyle{remark}
\newtheorem{rem}[thm]{Remark}
\newtheorem{ex}[thm]{Example}
\newtheorem{exs}[thm]{Examples}
\long\def\Thm#1{\begin{thm} #1 \end{thm}}
\long\def\Cor#1{\begin{cor} #1 \end{cor}}

\long\def\Lem#1{\begin{lem} #1 \end{lem}}
\long\def\Prop#1{\begin{prop} #1 \end{prop}}

\long\def\Def#1{\begin{defn} #1 \end{defn}}
\long\def\Rem#1{\begin{rem} #1 \end{rem}}
\long\def\Ex#1{\begin{ex} #1 \end{ex}}

\def\bar#1{\overline{#1}}
\def\Sect{\section}

\long\def\Ref#1#2#3#4#5#6{
\bibitem{#1}
{\rm #2,}
\textit{#3.}
{\rm #4}
\textbf{#5}
{\rm #6.}
}
\long\def\Refb#1#2#3#4{
\bibitem{#1}
{\rm #2,}
\textit{#3.}
#4.
}
\def\Zz{{\mathbb Z}}%blackboard bold Z
\def\Rr{{\mathbb R}}%blackboard bold R
\def\Cc{{\mathbb C}}%blackboard bold C
\def\Hh{{\mathbb H}}%blackboard bold H
\def\Ff{{\mathbb F}}

\def\SS{\mathfrak{S}}
\def\i{{\rm i}}
\def\e{{\rm e}}
\def\phi{\varphi}
\def\into{\hookrightarrow}
\def\iso{\cong}%isomorphism
\def\leq{\leqslant}
\def\geq{\geqslant}
\def\st{\mid}
\def\O{{\rm O}}
\def\U{{\rm U}}
\def\SO{{\rm SO}}
\def\PO{{\rm PO}}
\def\Sp{{\rm Sp}}
\def\PSO{{\rm PSO}}
\def\Hom{{\rm Hom}}
\def\comp{\circ}
\begin{document}

\title{On Borsuk--Ulam theorems and convex sets}

\author{M.~C.~Crabb}
\address{%
Institute of Mathematics\\
University of Aberdeen \\
Aberdeen AB24 3UE \\
UK}
\email{m.crabb@abdn.ac.uk}
\date{August 2022}
\begin{abstract}
The Intermediate Value Theorem is used to give an elementary
proof of a Borsuk-Ulam theorem of Adams, Bush and Frick \cite{ABF}
that, if $f: S^1\to \Rr^{2k+1}$ is a continuous function on
the unit circle $S^1$ in $\Cc$ such that
$f(-z)=-f(z)$ for all $z\in S^1$, then there is a finite
subset $X$ of $S^1$ of diameter at most $\pi -\pi /(2k+1)$
(in the standard metric in which the circle has
circumference of length $2\pi$)
such the convex hull of $f(X)$ contains $0\in\Rr^{2k+1}$.
\end{abstract}
 %characteristic classes
\maketitle
\Sect{Introduction}
We shall use the Intermediate Value Theorem to give an
elementary proof of the following Borsuk-Ulam theorem
of Adams, Bush and Frick in which $k\geq 1$ is a natural number,
$\zeta =\e^{2\pi\i/(2k+1)}\in\Cc$ 
and the metric $d$ on the unit circle,
$S(\Cc )$, in $\Cc$ is given by $d(z,\e^{\i\theta} z)=|\theta|$
if $|\theta | \leq \pi$.
\Thm{\label{main}
{\rm (\cite[Theorems 1 and 5]{ABF})}.
Let $f: S^1=S(\Cc )\to \Rr^{2k+1}$ be a continuous map such that
$f(-z)=-f(z)$ for all $z\in S(\Cc )$.
Then there exist $e_i\in\{ \pm 1\}$, $i=0,\ldots ,2k$,
and $z\in S(\Cc )$ such that $0$ lies in the convex hull
of the image $f(X)$ of the finite
set $X=\{ e_i\zeta^iz\st i=0,\ldots , 2k\}$,
which is a subset of $S(\Cc )$ with
diameter at most $\pi -\pi /(2k+1)$.

Moreover, there is an example of
such a map $f$ with the property that any
finite subset $X\subseteq S(\Cc )$ such that $0$ lies in the
convex hull of $f(X)$ has diameter greater than or equal to
$\pi -\pi /(2k+1)$.
}
Cohomological methods, using little more than knowledge of the
$\Ff_2$-cohomology of a real projective space, will also, 
in Section 3, establish a Borsuk--Ulam theorem for higher dimensional spheres strengthening \cite[Theorem 3]{ABF}.
The unit sphere $S(\Rr^n)$ in $\Rr^n$ is equipped with the 
standard metric
$d$: $\cos (d(u,v)) =\langle u,v\rangle$, $0\leq d(u,v)\leq\pi$.
\Thm{\label{second}
Let $m,\, n\geq 1$ be positive integers such that 
$m\leq 2^r\leq n<2^{r+1}$, where $r\geq 0$ is a non-negative integer.
Suppose that $f : S^{n-1}=S(\Rr^n)\to \Rr^{m+n-1}$ is a 
continuous map
such that $f(-v)=-f(v)$ for all $v\in S(\Rr^n)$. Then there
exists a finite subset $X\subseteq S(\Rr^n)$ with cardinality
at most $m+n$ and diameter at most $\pi -\arccos (1/n)$
such that $0$ lies in the convex hull of $f(X)$ in $\Rr^{m+n-1}$.

Moreover, there is an example of
such a map $f$ with the property that any
finite subset $X\subseteq S(\Rr^n )$ such that $0$ lies in the
convex hull of $f(X)$ has diameter greater than or equal to
$\pi -\arccos (1/n)$.
}
Consider, in general, a continuous map 
$f: S(\Rr^n) \to \Rr^{m+n-1}$, where $m\geq 0$ and $n >1$,
such that $f(-v)=-f(v)$ for all $v\in S(\Rr^n)$.
The classical Borsuk-Ulam theorem asserts,
in one form, that, if $m=0$,
there is a point $x\in S(\Rr^n)$ such that $f(x)=0$.
If $n=2$, the proof is an elementary exercise using
the Intermediate Value Theorem, while for general $n$
a proof (one of many) can be given using the
$\Ff_2$-cohomology of real projective space.
If $m >0$, the theorem clearly fails; for the inclusion
$S(\Rr^n)\subseteq \Rr^n\into \Rr^{m+n-1}$ has no zero.
But the Adams-Bush-Frick theorems show that there
is a finite subset $X\subseteq S(\Rr^n)$ with the property that zero
is a convex linear combination of the values $f(x)$,
$x\in X$, and, to make the assertion non-trivial (because
$\frac{1}{2}f(x)+\frac{1}{2}f(-x)=0$ for any $x$),
satisfying the condition that
$X$ does not contain any pair of antipodal points. 
And this condition is refined, for given $m$ and $n$, 
by bounding the diameter of $X$. A precise statement is given in 
Corollary \ref{one}.

For the wider context of these results the reader is referred
to \cite{ABF, ABF2}.

\medskip

{\it Acknowledgment.}  
I am grateful to the authors of \cite{ABF},
H.~Adams, J.~Bush and F.~Frick, for their helpful comments on
the 2019 version of this note.
\Sect{The Intermediate Value Theorem}
\Lem{\label{borsuk}
Let $f: S(\Cc )\to \Rr^{2k+1}$ be a continuous map
such that $f(-z)=-f(z)$.
Suppose that
$w_0,\ldots , w_{2k}$ are any $2k+1$ points in $S(\Cc)$.
Then there exist $e_i\in\{ \pm 1\}$, $\lambda_i\geq 0$,
for $i=0,\ldots , 2k$, with $\sum\lambda_i=1$,
and $z\in S(\Cc)$,
such that $\sum_{i=0}^{2k} \lambda_i f(e_izw_i)=0$.

If $f(w_0),\ldots , f(w_{2k})$ lie in a
$2k$-dimensional subspace of $\Rr^{2k+1}$, we can
require that $z=1$.
}
\begin{proof}
Consider the determinant map 
$$
\phi : S(\Cc ) \to \Lambda^{2k+1} \Rr^{2k+1},
z\mapsto f(zw_0)\wedge \cdots \wedge f(zw_{2k}).
$$
Then $\phi (-z) = (-1)^{2k+1}\phi (z)=-\phi (z)$.
So, by the Intermediate Value Theorem, $\phi$ has
a zero. 

If $\phi (z)=0$,
the vectors $f(zw_i)$ in $\Rr^{2k+1}$ are
linearly dependent and there exist $\mu_i\in\Rr$, not
all zero such that $\sum_i \mu_i f(w_i)=0$.
We may assume, by scaling, that $\sum |\mu_i|=1$.
Choose $\lambda_i\geq 0$ and $e_i=\pm 1$ so that
$\lambda_ie_i=\mu_i$.
\end{proof}
\Lem{\label{polygon}
For $z\in S(\Cc )$ and $e_i=\pm 1$, $i=0,\ldots ,2k$,
the distance $d(e_iz\zeta^i,e_jz\zeta^j)$ is
less than or equal to $\pi -\pi /(2k+1)$.
}
\begin{proof}
Indeed, the $2(2k+1)$ points $\pm z\zeta^i$ on the unit circle lie at
the vertices of a regular polygon.
\end{proof}
These two lemmas already prove the first part of Theorem \ref{main}.
\Ex{\label{interpolation}
(\cite[Theorem 5]{ABF}).
Let $P$ be the $k$-dimensional complex vector space of
complex polynomials $p(z)$ of degree $\leq 2k-1$
such that $p(-z)=-p(z)$. Let $g: S(\Cc ) \to P^*$ be
the evaluation map to the dual $P^*=\Hom_\Cc (P,\Cc )$.
Suppose that $w_0,\ldots , w_{2k}$ are points of $S(\Cc )$ such
that $0$ lies in the convex hull of the $g(w_i)$.
Then $d(w_i,w_j)\geq \pi -\pi /(2k+1)$ for some $i,\, j$. 
}
For the sake of completeness we include a concise version
of the proof in \cite{ABF}.
\begin{proof}
Assume that the $2k+1$ points $w_i^2$ are distinct.
By relabelling we may arrange that $w_1=\e^{\i\theta} w_0$,
where $\theta$ is the minimum of the distances
$d(w_i,w_j)$, $i\not=j$.

Suppose that $\sum\lambda_ig(w_i)=0$, where the $\lambda_i\in\Rr$ are
not all equal to zero. This means that $\sum \lambda_ip(w_i)=0$ for each $p\in P$, and hence, since $\lambda_i\in\Rr$ and 
$w_i^{-1}=\bar{w_i}$, that
$\sum \lambda_i p(w_i^{-1})=0$ too (because 
$\sum \lambda_i p(\bar{w_i})$ can be written as
the complex conjugate of $\sum \lambda_i \bar{p} (w_i)$
with $\bar{p}\in P$).
For $r\not=s$, $0\leq r,s\leq 2k$, we may write
$$
z^{-2k+1}\prod_{j\not=r,\, j\not=s} (z^2-w_j^2)=p_+(z)+p_-(z^{-1}),
$$
for unique polynomials $p_+,\, p_-\in P$.
Then we find, because $\sum \lambda_i p_+(w_i)=0$
and $\sum \lambda_i p_-(w_i^{-1})=0$, that
$$
\sum_i \lambda_i w_i^{-2k+1}\prod_{j\not=r,\, j\not=s} (w_i^2-w_j^2)=
\sum_i \lambda_i (p_+(w_i)+ p_-(w_i^{-1}))=0,
$$
that is,
$$
\lambda_r\prod_{j\not=r} (w_rw_j^{-1}-w_r^{-1}w_j)
=\lambda_s\prod_{j\not=s} (w_sw_j^{-1}-w_s^{-1}w_j)\, .
$$
It follows that, for some non-zero $c\in\Rr$,
$$
\lambda_i \delta_i =c,
\text{\ where\ }\delta_i=\prod_{j\not=i} (w_iw_j^{-1}-w_i^{-1}w_j),
$$
for all $i$. (Notice that $\delta_i$, being the product of the
$2k$ purely imaginary numbers $w_i\bar{w_j}-\bar{w_i}w_j$, is real.)
In particular, all the $\lambda_i$ are non-zero.

Given that $0$ lies in the convex hull of the points $g(w_i)$,
we can now assume further that all $\lambda_i$ are non-negative,
and so, because they are non-zero, strictly positive. Then the
$\delta_i \in\Rr$ all have the same sign.
We show that there is some $i$ such that 
$w_i=-\e^{\i t\theta}w_0$ for $0<t<1$.

Indeed, write $\psi (t)= \prod_{1<j\leq 2k} 
(\e^{\i t\theta}w_0w_j^{-1}
-\e^{-\i t\theta}w_0^{-1}w_j)\in\i\Rr$, for $0\leq t\leq 1$.
Then $\delta_0 = (w_0w_1^{-1}-w_0^{-1}w_1)\psi (0)$
and $\delta_1= (w_1w_0^{-1}-w_1^{-1}w_0)\psi (1)$.
So, by the Intermediate Value theorem again,
$\psi (t)=0$ for some $t$, and then 
$w_i^2 =(\e^{\i t\theta}w_0)^2$ for some $i$, $1<i\leq 2k$.
But $w_i\not= \e^{\i t\theta}w_0$, by the minimality of
$\theta$. So $w_i=-\e^{\i t\theta}w_0$.
Now $d(w_i,w_0)=\pi -t\theta$ and
$d(w_i,w_1)=\pi -(1-t)\theta$. But clearly $\theta
\leq 2\pi /(2k+1)$ and either $t\geq 1/2$ or
$1-t\geq 1/2$.
\end{proof}
Lemmas \ref{borsuk} and \ref{polygon}, together with Example
\ref{interpolation}, establish Theorem \ref{main}, using the 
fact\footnote{We have $\lambda_y>0$, $y\in Y$, such  that $\sum \lambda_y y=0$. Suppose that $\mu_y\in\Rr$, $y\in Y$, satisfy
$\sum\mu_y y=0$ and $\sum \mu_y=0$. Then,
for any $t\in\Rr$ such that $\lambda_y\geq |t\mu_y|$
for all $y$, $\sum (\lambda_y+t\mu_y)y=0$, 
$\sum (\lambda_y-t\mu_y)y=0$, $\lambda_y\pm t\mu_y\geq 0$,
and so $\lambda_y >\pm t\mu_y$, that is, $\lambda_y>|t\mu_y|$,
for all $y$. 
Hence $\mu_y=0$ for all $y$.}
that, if $Y$ is a finite subset of a real vector space such  
that the convex hull of $Y$ contains $0$ but no proper  
subset of $Y$ has this property, then the set $Y$ is 
affinely independent.
\Cor{\label{one}
Let $m\geq 0$ and $n >1$ be integers. Then there is a non-negative
real number
$\delta <1$ with the property that, for any continuous map 
$f: S(\Rr^n)\to\Rr^{m+n-1}$ such that $f(-v)=-f(v)$ for all
$v\in S(\Rr^n)$, there is a finite subset $X$ of $S(\Rr^n)$
of cardinality at most $m+n$ and diameter at most 
$\pi -\arccos(\delta )<\pi $ such that $0$ lies in the convex hull of 
$f(X)$.
}
The classical Borsuk-Ulam theorem deals with the case
$m=0$ : we may take $\delta =0$ so that
$X$ consists of a single point.
\begin{proof}
By restricting $f$ to $S(\Rr^2)\subseteq S(\Rr^n)$ and including 
$\Rr^{m+n-1}$ in $\Rr^{2k+1}$ for the smallest $k\geq 1$ such that
$m+n-1\leq 2k+1$, we see from Theorem \ref{main} that the
assertion is true with $\delta =\cos (\pi /(2k+1))$.
\end{proof}
There is an easy extension of Lemma \ref{borsuk} to higher dimensions.
\Lem{\label{gen1}
For integers $n,\, k\geq 1$, 
suppose that $f : S(\Cc^{n})\to
\Rr^{2k+1}$ is a continuous map such that $f(-v)=-f(v)$
for all $v\in S(\Cc^{n})$.
Then, for any $2k+1$ vectors $w_0,\ldots ,w_{2k}$
in $S(\Cc^{n})$, there exist $e_i\in\{\pm 1\}$, $\lambda_i\geq 0$, 
for $i=0,\ldots 2k$, with $\sum\lambda_i=1$, and
$z\in S(\Cc)$ such that $\sum_{i=0}^{2k} \lambda_i f(e_izw_i)=0$.
}
\begin{proof}
This can be established, using the same arguments as in the proof of 
Lemma \ref{borsuk}, by looking at the function 
$\phi : S(\Cc ) \to \Lambda^{2k+1} \Rr^{2k+1}$ defined by 
$\phi (z)= f(zw_0)\wedge \cdots\wedge f(zw_{2k})$. 
\end{proof}
\Ex{\label{mult}
For $n >1$, write ${\bf e}_1,\ldots ,{\bf e}_n$ for the standard 
orthonormal $\Cc$-basis of $\Cc^n$,
let $\eta =\e^{2\pi\i /(2l+1)}$ where $l\geq 1$ is a positive
integer, and fix an integer $r$ in
the range $1\leq r\leq n$.

If $k$ satisfies $2k+1\leq (2l+1)^r\binom{n}{r}$,
we can choose distinct vectors $w_0.\ldots ,w_{2k}$ from the set
$$\textstyle
\{(\sum_{s=1}^r\eta^{a_s}{\bf e}_{i_s})/\sqrt{r}\st 
1\leq i_1<i_2\cdots < i_r\leq n,\,
a_s=0,\ldots ,2l\}\subseteq S(\Cc^n).
$$
Then $|\langle w_i,w_j\rangle | \leq 
\delta =1-(1-\cos (\pi /(2l+1))/r$ for $i\not=j$. 
Hence $d(e_izw_i,e_jzw_j)\leq 
\pi -\arccos (\delta )$.

Application of Lemma \ref{gen1} to the case
$r=1$, for which $\arccos (\delta )=\pi /(2l+1)$,
gives a result which is close to
\cite[Theorem 2]{ABF}, but slightly weaker.
For general $r$, if $l$ is large, $\arccos (\delta )$ is close to
$(\pi /(2l+1))/\sqrt{r}$. In particular, the case
$r=n$ gives a much stronger result than that provided by
\cite[Theorem 2]{ABF} when $k$ is sufficiently large.
(A similar observation is made in \cite[Remark 4.2]{ABF2}.)
}
\Sect{Cohomological methods}
In this section we shall use cohomological methods in higher 
dimensions. We begin by recalling some standard facts about
projective bundles. 

Consider a real vector bundle $E\to B$ (or simply `$E$')
of dimension $n$ over a 
finite complex $B$. The {\it mod $2$ Euler class} of $E$ will
be written as 
$$
e(E)\in H^n(B;\,\Ff_2)
$$
in cohomology with $\Ff_2$-coefficients. (Thus, $e(E)$ is the top
Stiefel-Whitney class $w_n(E)$.) The crucial property for
our applications is that, if $e(E)\not=0$,
then every section $s:B\to E$ has a zero, that is, there is 
some point $b\in B$ such that $s(b)\in E_b$, the vector space
fibre of $E$ at $b$, is equal to $0\in E_b$.
For calculations we shall use the {\it multiplicativity}
of the Euler class: $e(E\oplus E')=e(E)\cdot e(E')$ for
two vector bundles $E$ and $E'$ over $B$.

The {\it projective bundle} $P(E)\to B$ of $E$ is the quotient
of the sphere bundle $S(E)$ of $E$ by the involution 
$-1$: $P(E)=S(E)/\{ \pm 1\}$;
elements of $P(E)$ can be written as $[u]$, where $u\in S(E)$.
The {\it Hopf line bundle} over $P(E)$ is denoted by $H$;
its fibre at $[u]$ is the line $\Rr u$ generated by $u$.

The cohomology of the space $P(E)$, as an algebra over
$H^*(B;\,\Ff_2)$, is a quotient of the polynomial
ring $H^*(B;\, \Ff_2)[T]$ in an indeterminate $T$:
$$
H^*(P(E);\, \Ff_2)=H^*(B;\, \Ff_2)[T]/(T^n+\ldots +
w_i(E)T^{n-i}+\ldots +w_n(E)),
$$
where $w_i(E)\in H^i(B;\,\Ff_2)$ is the $i$th 
{\it Stiefel-Whitney class} of $E$, and
the Euler class $e(H)$ corresponds to $T$. 

The multiplicativity of the Euler class translates into
the multiplicativity of the {\it total Stiefel-Whitney class}
$w(E)= 1+w_1(E)+\ldots +w_i(E)+\ldots +w_n(E)$:
$w(E\oplus E')=w(E)\cdot w(E')$.

Working now towards a proof of Theorem \ref{second} we consider,
for $n>1$, an inscribed regular $n$-simplex in the sphere
$S(\Rr^n)$ with vertices $v_0,\ldots ,v_{n}$.
Write $n=2^r+s$, where $r\geq 1$ and $0\leq s<2^r$.
An involution $\tau$ in the orthogonal
group $\O (\Rr^n)$ is specified by: $\tau (v_i)=v_{i+s}$ for
$i<s$ and $\tau (v_i)=v_i$ for $i=2s,\ldots ,n$.
Let $V_+$ and $V_-$ be the $\pm 1$-eigenspaces of $\tau$:
$V_{\pm}=\{ u\in\Rr^n\st \tau u=\pm u\}$. Thus $V_+$ has dimension
$n-s=2^r$ with a basis $v_i+v_{i+s} \, (0\leq i<s)$,
$v_i\, (2s < i\leq n)$, and $V_-$ has dimension $s$ with
a basis $v_i-v_{i+s}\, (0\leq i<s)$. 
The special orthogonal group $\SO (V_+)$,
which is isomorphic to $\SO (\Rr^{2^r})$, is included
in $\SO (\Rr^n)$ by the map
$g\mapsto (g,1): V_+\oplus V_-=\Rr^n\to V_+\oplus V_-=\Rr^n$
as the subgroup $G\leq \SO (\Rr^n)$ 
consisting of those elements which fix
the $s$ vectors $v_i-v_{i+s}$, $0\leq i<s$.

\Thm{\label{simplex}
{\rm (Compare \cite[Theorem 3]{ABF})}.
For an integer $n > 1$,
let $v_0,\ldots ,v_{n}$ be the vertices of an inscribed regular 
$n$-simplex in $S(\Rr^{n})$.
Write $n=2^r+s$, where $r\geq 1$ and $0\leq s<2^r$,
and let $G\leq \SO (\Rr^n)$ be the subgroup 
$\{ g\in\SO (\Rr^n)\st g(v_i-v_{i+s})=v_i-v_{i+s}\, 
\text{for $0\leq i<s$}\}$.

Let $f: S(\Rr^{n})\to \Rr^{n+2^{r}-1}$ 
be a continuous map such that
$f(-v)=-f(v)$ for all $v\in S(\Rr^{n} )$.
Then there exist $e_i\in\{ \pm 1\}$, $i=0,\ldots ,n$,
and an element $g\in G$ such that $0$ lies in the convex hull
of the image $f(X)$ of the finite
subset $X=\{ e_i gv_i\st i=0,\ldots , n\}\subseteq S(\Rr^{n} )$,
which has diameter at most $\pi -\arccos (1/n)$.

Moreover, the example of the inclusion
$f: S(\Rr^{n}) \subseteq \Rr^{n} \into \Rr^{n+2^r-1}$
has the property that any
finite subset $X\subseteq S(\Rr^n )$ such that $0$ lies in the
convex hull of $f(X)$ has diameter greater than or equal to
$\pi -\arccos (1/n)$.
}
\begin{proof}
We continue to use the notation introduced before the
statement of the theorem.

Consider the map $\sigma : \SO (V_+)\times S(\Rr^{n+1})
\to \Rr^{n+2^{r}-1}$
$$
(g,(\mu_0,\ldots ,\mu_{n}))\mapsto
\sum_{i=0}^{n} \mu_i f(gv_i), \quad g\in\SO (V_+),\quad
\sum_{i=0}^n \mu_i^2=1.
$$
We shall show that $\sigma$ has a zero.

Clearly $\sigma (g,-\mu )=-\sigma (g, \mu )$.
Now notice that $-g\in \SO (V_+)$ maps to
$(-g,1) = -(g,1)\tau\in \SO(\Rr^n)$. So
$\sigma (-g,\mu ) = -\sigma (g,T\mu)$, where
$T$ is the linear involution of $\Rr^{n+1}$ satisfying,
for $\mu=(\mu_0,\,\ldots ,\, \mu_n)\in\Rr^{n+1}$,
$(T\mu)_i =\mu_{i+s}$ for $i<s$
and $(T\mu)_i =\mu_i$ for $i=2s,\ldots ,n$.

Let $L$ be the real line bundle over the projective special
orthogonal group $\PSO (V_+)$ associated with the double cover
$\SO (V_+)\to \PSO (V_+)=\SO (V_+)/\{ \pm 1\}$ or, explicitly,
the quotient of $\SO (V_+)\times \Rr$ by the involution
$(g,t)\mapsto (-g.-t)$.
And let $E$ be the $(n+1)$-dimensional real vector bundle
over $\PSO(V_+)$ associated to the involution $-T$ of $\Rr^{n+1}$,
that is, $E$ is the quotient of $\SO (V_+)\times \Rr^{n+1}$
by the involution $(g,\mu )\mapsto (-g,-T\mu )$.
The $-1$-eigenspace of $-T$ consists of those $\mu\in\Rr^{n+1}$
such that $\mu_i=\mu_{i+s}$ for $0\leq i<s$ and has 
dimension $n+1-s=2^r+1$, and the $+1$-eigenspace has 
the complementary dimension $s$.
The bundle $E$ is thus isomorphic to the direct sum
$(L\otimes\Rr^{2^r+1})\oplus \Rr^s$ of $2^r+1$ copies of $L$
and the trivial bundle with fibre $\Rr^s$.

The map $\sigma$ determines a section $s$ of $H\otimes\Rr^{n+2^r-1}$
(the tensor product of the
pullback of the Hopf line bundle $H$ over $P(\Rr^{n+1})$
with the trivial bundle with fibre $\Rr^{n+2^r-1}$):
$$
s([[g,\pm\mu ]]) = [g,\mu]\otimes\sigma (g,\mu ),\quad
\text{for $g\in\SO (V_+)$, $\mu\in S(\Rr^{n+1})$,}
$$
where the element $[[g,\mu ]]=[[g,-\mu ]]\in P(E)$ is 
determined by the class 
$[g,\mu ]\in S(E)$ of $(g,\mu )\in \SO (V_+)\times S(\Rr^{n+1})$.
We shall show that $s$ has a zero, and hence that $\sigma$ has a 
zero, by proving that the mod $2$ Euler class 
$$
e(H\otimes\Rr^{n+2^r-1})=e(H)^{n+2^r-1}\in
H^{n+2^r-1}(P(E);\,\Ff_2)
$$
is non-zero.

Let us write $x=e(L)$ and $y=e(H)$. 
Now we recall from \cite[Lemma (8.1)]{proj} that
$x^{2^{r}-1}\not=0$ and $x^{2^r}=0$. 
(In order to make this account as self-contained as possible,
an elementary proof of this fact is included in an Appendix as 
Proposition \ref{PSO}.)
So the total Stiefel-Whitney class $w(E)$ of $E$ is determined,
using multiplicativity, by $w(L)=1+x$ as
$w(E)=w(L)^{2^r+1}=(1+x)^{2^r+1}=1+x$ 
(because $(1+x)^{2^r}=1+x^{2^r}=1$). Hence
$$
H^*(P(E);\,\Ff_2) =H^*(\PSO (V_+);\, \Ff_2)[y]/(y^{n+1}+xy^n),
$$
(meaning that the quotient
of the polynomial ring $H^*(\PSO (V_+);\,\Ff_2)[T]$ by the principal 
ideal generated by $T^{n+1}+xT^n$ is isomorphic to
$H^*(P(E);\,\Ff_2)$ by the homomorphism mapping 
$T$ to $y$).
Using the relation $y^{n+1}=xy^n$,
we conclude that $y^{n+2^r-1}=x^{2^r-1}y^n$ is non-zero.

Hence $\sigma$ has a zero, that is, there exist $g\in\SO (V_+)$
and $\mu_0,\,\ldots ,\,\mu_n$, not all zero, such that 
$\sum_i \mu_if(gv_i)=0$. The existence of the required
$e_i\in\{\pm 1\}$ follows as in the proof of Lemma \ref{borsuk}.

The $n+1$ points $v_i$ satisfy $\sum_i v_i=0$ and 
$\langle v_i,v_j\rangle =c$, say, for $i\not=j$.
Since $\sum_{i=0}^n \langle v_0,v_i\rangle =0$, we see that
$1+nc=0$, so that $c=-1/n$.
Hence $d(v_i,v_j)=\pi -\arccos (1/n)$ and
$d(v_i,-v_j)=\arccos (1/n)$ if $i\not=j$.

So $d(e_igv_i,e_jgv_j)=d(e_iv_i,e_jv_j)\leq \pi -\arccos (1/n)$
for $i\not=j$.

The final assertion follows, as in the proof of Theorem \ref{main},
from Lemma \ref{below} below.
\end{proof}
\Lem{\label{below}
Suppose that $w_0,\ldots , w_{n}$ are points of 
$S(\Rr^n)$, such that $\sum \lambda_i w_i=0$, where $\lambda_i\geq 0$
and $\sum\lambda_i=1$. Then 
$d(w_i,w_j)\geq \pi -\arccos (1/n)$ for some $i,\, j$.
}
\begin{proof}
Put $d=\min \{ \langle w_i,w_j\rangle \st
i\not=j\}$. Then
$$
0=\sum_i \sum_j \lambda_j\langle w_i,w_j\rangle
=\sum_i\lambda_i\langle w_i,w_i\rangle +
\sum_i\sum_{j\not=i} \lambda_j\langle w_i,w_j\rangle \geq
1+ nd\, .
$$
So $d\leq -1/n$. If $\langle w_i,w_j\rangle\leq -1/n$, then
$d(w_i,w_j)\geq \pi -\arccos (1/n)$.
\end{proof}
The result in Theorem \ref{simplex} is 
optimal when $n=2^r$ is a power of $2$ and so $n+2^r-1=2n-1$,
as the following construction demonstrates.
\Prop{\label{even}
For any even integer $n > 1$, there is a map
$f :S(\Rr^n) \to \Rr^{2n}$ such that for every
finite subset $X$ of $S(\Rr^n)$ with diameter
$\leq \pi -\arccos (1/n)$ the zero vector
$0\in\Rr^{2n}$ does not lie in the convex hull of $f(X)$.
}
\begin{proof}
Consider a map $f : S(\Rr^n) \to \Rr^{2n}=\Rr^n\oplus\Rr^n$
of the form $f(v)=(v,\phi (v)v)$, where $\phi : S(\Rr^n)\to
\Rr$ is a continuous function such that $\phi (-v)=\phi (v)$.

Suppose that $X\subseteq S(\Rr^n)$ is a finite subset of
diameter at most $\pi -\arccos (1/n)$ such that
$0=(0,0)$ lies in the convex hull of $f(X)$.

Considering the first component, we see that there
are points $w_0,\ldots , w_m$ of 
$X$ for some $m$ with $1\leq m\leq n$ 
and $\lambda_0,\ldots ,\lambda_m\geq 0$
with $\sum\lambda_i=1$, such that $\sum \lambda_i w_i=0$ and 
$\langle w_i,w_j\rangle \geq -1/n$ for all $i,\, j$.  
(Indeed, let $Y\subseteq X$ be a minimal subset of $X$
such that $0$ lies in the convex hull of $Y$.
As we have already recalled in the proof of Theorem \ref{main},
the set $Y$ must be affinely independent and so has cardinality,
$m+1$ say, less than or equal to $n+1$.)

Then, for each $i$,
$$
-\lambda_i=\sum_{j\not=i}\lambda_j \langle w_i,w_j\rangle 
\geq -(1-\lambda_i)/n,
$$
and so $\lambda_i\leq 1/(n+1)$.
Hence $m=n$ and $\lambda_i=1/(n+1)$ for all $i$. And then
$\sum_{j\not=i} \langle w_i,w_j\rangle =-1$.
We conclude that $\langle w_i,w_j\rangle =-1/n$
for $i\not=j$. Hence $w_0,\ldots ,w_n$ are the vertices
of an inscribed regular $n$-simplex in $S(\Rr^n)$.

If $v\in S(\Rr^n)$ is a vector distinct from each of
$w_0,\ldots ,w_n$, then $\langle v,w_i\rangle <-1/n$
for some $i$. 
(Indeed, if $v=\sum_{j\not=i}t_jw_j$,
where $t_j\geq 0$, then $\langle v,w_i\rangle
< -1/n$.
For $1=\| v\|^2 
=(\sum_{j\not=i}t_j)^2-
((n+1)/n)\sum_{j\not=j',\, j\not=i\not=j}t_jt_{j'}$.
So $\sum_{j\not=i}t_j >1$. Hence
$\langle v,w_i\rangle = -(\sum_{j\not=i}t_j)/n < -1/n$.)
Thus $X$ is precisely the set $\{ w_0,\ldots ,w_n\}$.

Since the $w_i$ are the vertices of a regular $n$-simplex,
if $\sum_i\mu_iw_i=0$, $\mu_i\in\Rr$, then the coefficients
$\mu_i$ must be equal.
So, looking now at the second component of $f$, we conclude
that the real numbers $\phi (w_i)$, $i=0,\ldots , n$, are
equal.
(For  $\sum_i \lambda_i f(w_i)=0$ means that $\sum_i \lambda_i w_i=0$
and $\sum_i \lambda_i \phi (w_i)w_i=0$, and 
we have shown that $\lambda_i=1/(n+1)$.)

To construct a map $f$ such that $0$ does not lie in the convex
hull of $f(X)$ for any finite set with diameter at most
$\pi -\arccos (1/n)$, 
we write down a function
$\phi$ which is not constant on the vertices $w_0,\,\ldots,\, w_n$ of
any regular $n$-simplex. For example, choose a non-zero
vector $u\in \Rr^n$ and take
$\phi (v) =\langle u,v\rangle^2$.
Then $0=\sum \langle u,w_i\rangle
= \sum_{i=0}^n \epsilon_i\sqrt{\phi (w_i)}$,
where $\epsilon_i=\pm 1$. Since
$n+1$ is odd, if $\phi (w_i)=c^2$ for all $i$, then
$c=0$, which is impossible because $u\not=0$.
\end{proof}
\begin{proof}[Proof of Theorem \ref{second}]
This follows at once from Theorem \ref{simplex} by including $\Rr^{m+n-1}$
in $\Rr^{n+2^r-1}$.
\end{proof}
Methods similar to those used to derive Theorem \ref{simplex}
give the following generalization of Lemma \ref{gen1}.
\Lem{\label{gen}
For some integers $r\geq 0$ and $n,\, k\geq 1$, 
suppose that $f : S(\Rr^{2^rn})\to
\Rr^{2^rk+2^r-1}$ is a continuous map such that $f(-v)=-f(v)$
for all $v\in S(\Rr^{2^rn})$.
We include the special orthogonal group
$\SO (\Rr^{2^r})$ in $\SO (\Rr^{2^rn})$ by the diagonal
homomorphism
$\SO (\Rr^{2^r}) \to \SO (\Rr^{2^r})\times\cdots\times 
\SO (\Rr^{2^r})\subseteq \SO (\Rr^{2^rn})$.
Then, for any $2^rk+1$ vectors $w_0,\ldots ,w_{2^rk}$
in $S(\Rr^{2^rn})$, there exist $e_i\in\{\pm 1\}$, $\lambda_i\geq 0$, 
for $i=0,\ldots 2^rk$, with $\sum\lambda_i=1$, and
$g\in \SO (\Rr^{2^r})$, such that $\sum_{i=0}^{2^rk} 
\lambda_i f(e_igw_i)=0$.

If, for some $\delta <1$, the vectors $w_i$ satisfy
$|\langle w_i,w_j\rangle | \leq \delta$ for all $i\not=j$,
then $d(e_igw_i,e_jgw_j)\leq \pi -\arccos (\delta )<\pi$ for $i\not=j$.
}
\begin{proof}
We assume that $r\geq 1$.
The simpler case $r=0$ is discussed in Remark \ref{r0}. 
Consider the map $\sigma : \SO (\Rr^{2^r})\times S(\Rr^{2^rk+1})
\to \Rr^{2^rk+2^{r}-1}$
$$
(g,(\mu_0,\ldots ,\mu_{2^rk}))\mapsto
\sum_{i=0}^{2^rk} \mu_i f(gw_i), \quad g\in\SO (\Rr^{2^r}),\quad
\sum_{i=0}^{2^rk} \mu_i^2=1.
$$
Since $\sigma (-g,\mu )=-\sigma (g, \mu ) = \sigma (g,-\mu )$,
the map $\sigma$ determines a section $s$ of the vector bundle
$\Rr^{2^rk+2^{r}-1}\otimes (L\otimes H)$ over the space
$(\SO (\Rr^{2^r})/\{ \pm 1\} )\times (S(\Rr^{2^rk+1})/\{ \pm 1\})$,
where $L$ and $H$ are the line bundles associated with the
double covers 
$\SO (\Rr^{2^r}) \to \SO (\Rr^{2^r})/\{\pm 1\}$  
and $S(\Rr^{2^rk+1})\to S(\Rr^{2^rk+1})/\{ \pm 1\}$
of $\PSO (\Rr^{2^r})$ and the
real projective space $P(\Rr^{2^rk+1})$ , respectively.

The $\Ff_2$-Euler class of $\Rr^{2^rk+2^{r}-1}\otimes (L\otimes H)$ 
is equal to 
$$
e(L\otimes H)^{2^rk+2^r-1} =(x+y)^{2^rk+2^{r}-1}
\in H^{2^rk+2^r-1}(\PSO (\Rr^{2^r})\times P(\Rr^{2^rk+1});\,\Ff_2),
$$
where $x=e(L)$ and $y=e(H)$. 
Since $x^{2^{r}-1}\not=0$ and $x^{2^r}=0$ (by Proposition \ref{PSO}
again),
$y^{2^rk}\not=0$ and $y^{2^rk+1}=0$, 
$$
(x+y)^{2^rk+2^{r}-1}=
(x+y)^{2^r-1}(x^{2^r}+y^{2^r})^k =(x+y)^{2^r-1}y^{2^rk}
=x^{2^r-1}y^{2^rk}
$$
is non-zero. Hence the section $s$ has a zero, and, 
therefore, so too has $\sigma$. The proof is completed as before.
\end{proof}
\Rem{\label{r0}
When $r=0$, so that $g=1$ (because $\SO (\Rr^1)$ is trivial),
the result specializes to the main theorem of \cite[Theorem 1]{ABF2}.
In that case,
the map $\sigma$ determines a section $s$ of $\Rr^k\otimes H$ over
$P(\Rr^{k+1})$, and $s$ has a zero just because $y^k=e(H)^{k}$
is non-zero in $H^{k}(P(\Rr^{k+1});\,\Ff_2)$.

The case $r=1$ is precisely Lemma \ref{gen1} with 
the identifications $\Rr^{2n}=\Cc^n$ and 
$\SO (\Rr^2)=\U (\Cc )=S(\Cc )$.

When $r=2$, we can write $\Rr^{4n}=\Hh^n$ and use quaternion
multiplication to show that the element $g\in \SO (\Rr^4)$
may be chosen to lie in the subgroup $\Sp (\Hh )=S(\Hh )$
of unit quaternions. (And, when $r=3$, we could use Cayley
multiplication.)
}
\Cor{\label{bound}
Let $m,\, n\geq 1$ be integers with $m\leq n$.
Then there are $m+n$ vectors $w_0,\ldots , w_{m+n-1}\in S(\Rr^n)$
with the property that, for any continuous map $f:S(\Rr^n)\to 
\Rr^{m+n-1}$ with $f(-v)=-f(v)$, there exist $e_0,\ldots ,e_{m+n-1}\in\{ \pm 1\}$ such that  $0$ lies in the convex
hull of the image $f(X)$ of the finite set
$X=\{ e_iw_i\st i=0,\ldots , m+n-1\}$
with cardinality $m+n$ and diameter at most
$\pi -\arccos (1/\lfloor n/m\rfloor)$.
}
Notice that the vectors $w_0,\ldots ,w_{m+n-1}$ are independent
of $f$; compare Theorem \ref{second}.
The construction we use can be seen as
a special case of \cite[Lemma 16]{bukh}.
Compare, also, Example \ref{mult}.
\begin{proof}
Write $n=qm+l$, where $q=\lfloor n/m\rfloor\geq 1$ and $0\leq l<m$.
We can assume that $q>1$.
Let $u_0,\ldots , u_q$ be the vertices of an inscribed regular
$q$-simplex in $S(\Rr^q)$
and let ${\bf e}_1,\ldots , {\bf e}_m$ be an 
orthonormal basis of $\Rr^m$.
We give $\Rr^q\otimes\Rr^m$ the standard inner product:
$\langle u\otimes v,u'\otimes v'\rangle =\langle u,u'\rangle
\cdot \langle v,v'\rangle$.
In other words, $\Rr^q\otimes\Rr^m$ is the orthogonal direct sum
$\bigoplus_{j=1}^m \Rr^q\otimes\Rr {\bf e}_j$.
Thus $\langle u_i\otimes {\bf e}_j,u_{i'}\otimes 
{\bf e}_{j'}\rangle$ is equal
to $1$ if $i=i'$, $j=j'$, equal to $-1/q$ if $i\not=i'$,
$j=j'$, and is zero otherwise.

Let $w_0,\ldots , w_{k-l}\in \Rr^{qm}\subseteq 
\Rr^{qm}\oplus\Rr^l=\Rr^n$, where $k=n+m-1$,
correspond to the $(q+1)m=n+m-l=k-l+1$ vectors
$u_i\otimes {\bf e}_j$, in some order, under a chosen
isometric isomorphism $\Rr^{qm}\iso \Rr^q\otimes\Rr^m$
and let $w_{k-l+1},\ldots ,w_{k-l+l}$ be an orthonormal
basis of $\Rr^l\subseteq \Rr^{qm}\oplus\Rr^l =\Rr^n$.
The result follows from Lemma \ref{gen} with $r=0$
and $\delta =1/q$.
\end{proof} 
\Sect{Conclusion}
We finish with a reformulation of the main results.
\Def{Given integers $m\geq 1$ and $n>1$, set $k=m+n-1$.
Let $\delta (m,n)$ be the smallest real number $\delta \in (0,1]$
with the property that, for any continuous
map $f : S(\Rr^n)\to\Rr^{k}$ satisfying
$f(-v)=-f(v)$ for all $v\in S(\Rr^n)$,
there is a finite subset $X\subseteq S(\Rr^n)$ of
cardinality less than or equal to $m+n$ and diameter at most 
$\pi -\arccos (\delta )$
such that $0$ lies in the convex hull of $f(X)$.

The existence of $\delta (m,n)$ follows from a compactness argument.
For a given $f$, the subspace
$$
\{(\lambda_0,\ldots , \lambda_k, w_0,\ldots ,w_k)\in [0,1]^{k+1}\times
S(\Rr^n)^{k+1} \st \sum_{i=0}^k\lambda_i=1,\, 
\sum_{i=0}^k \lambda_i f(w_i)=0\}
$$
of $[0,1]^{k+1}\times S(\Rr^n)^{k+1}$ is compact and non-empty.
So the subset $\Delta (f)\subseteq [0,1]$ consisting of those 
$\delta$ for which
there is a subset $X$ of diameter $\leq \pi-\arccos(\delta )$
such that $0$ is in the convex hull of $f(X)$ is closed. Hence
$\bigcap_f\Delta (f)$ is closed (and contains $1$).

The example of the inclusion
of $S(\Rr^n)$ in $\Rr^n\subseteq\Rr^k$ shows, by Lemma \ref{below},
that $\delta (m,n)\geq 1/n$. Corollary \ref{one} asserts that
$\delta (m,n)<1$.
}
\Thm{\label{summary}
The function $\delta (m,n)$ defined for $m\geq 1$ and $n >1$ 
has the following properties.
\smallskip

\par\noindent
{\rm (i).} For all $m$ and $n$,
$1/n\leq \delta (m,n) <1$.

\smallskip

\par\noindent
{\rm (ii).}
For all $k\geq 1$,
$\delta (2k-1,2)=\delta (2k,2)=\cos (\pi /(2k+1))$
{\rm (Theorem \ref{main}, Example \ref{interpolation}).} 

\smallskip

\par\noindent
{\rm (iii).} 
If there is an integer $r$ such that $1\leq m\leq 2^r\leq n$,
then $\delta (m,n)=1/n$
{\rm (Theorem \ref{second}).}

\smallskip

\par\noindent
{\rm (iv).}
For all $k\geq 1$,
$\delta (2k+1,2k) >1/(2k)$
{\rm (Proposition \ref{even}).}

\smallskip

\par\noindent
{\rm (v).} 
Fix $m\geq 1$. Then
$\delta (m,n)\to 0$ as $n\to\infty$,
and so $\pi -\arccos (\delta (m,n))\to \pi /2$.
{\rm (Compare \cite[Corollary 3.1]{ABF2}.)}

\smallskip

\par\noindent
{\rm (vi).}
Fix $n>1$.
Then $\delta (m,n) \to 1$ as $m\to\infty$,
and so $\pi -\arccos(\delta (m,n)) \to \pi$.
}
\begin{proof}
Part (v) follows at once from (iii).
Only (vi) still requires proof.
If $m\leq m'$, then clearly $\delta (m,n)\leq \delta (m',n)$.

Given $n$ and a positive integer $l\geq 1$, set
$m+n-1=k=\binom{n+2l}{n-1}$.
We shall establish the easy estimate 
$\delta (m,n) \geq 1/k^{1/(2l+1)}$.
Since $k$, the dimension of the space of homogeneous polynomials
of degree $2l+1$ in $n$ variables, is clearly bounded by
$2l+1\leq k\leq (2l+1)^{n-1}$, we see that,
as $l\to\infty$ (with $n$ fixed), $m\to\infty$ and 
$k^{1/(2l+1)}\to 1$.

Consider the map $f : S(\Rr^n)\to (\Rr^n)^{\otimes (2l+1)}$ given
by the $(2l+1)$-th power $f(v)=v\otimes \cdots\otimes v$
and mapping into the $k$-dimensional vector subspace fixed
by the permutation action of $\SS_{2l+1}$.
Give $(\Rr^n)^{\otimes (2l+1)}$ the standard inner product:
$$
\langle u_1\otimes\cdots\otimes u_{2l+1},v_1\otimes\cdots
\otimes v_{2l+1}\rangle =\prod_i \langle u_i,v_i\rangle,
$$
so that each vector $f(v)$ lies in the unit sphere. 

By Lemma \ref{below}, if $w_0,\ldots ,w_k$ are points
of $S(\Rr^n)$ such that $0$ lies in the convex hull
of the unit vectors $f(w_0),\ldots ,f(w_k)$, then
$\langle f(w_i),f(w_j)\rangle = \langle w_i,w_j\rangle^{2l+1}
\leq -1/k$ for some $i\not=j$.
So $\langle w_i,w_j\rangle \leq -1/k^{1/(2l+1)}$. 
\end{proof}
\Rem{The example used to prove part (vi) of Theorem \ref{summary}
shows more precisely that, for any real number $\alpha<1/(n-1)$,
$$
m^\alpha (1-\delta (m,n)) \to 0
\text{\ as\ } m\to\infty\, .
$$
}
\begin{proof}
The construction shows that $\delta ((2l+1)^{n-1},n) \geq 
1/(2l+1)^{(n-1)/(2l+1)}$. Write $x=1-\delta ((2l+1)^{n-1},n)$. 
Then 
$$
(2l+1)^{n-1}\geq (1-x)^{-(2l+1)}=\sum_{s\geq 0} \binom{2l+s}{s}x^s
\geq \binom{2l+r}{r}x^r\geq (2l+1)^rx^r/r!
$$ 
for any $r\geq 1$.
Thus $x\leq (r!)^{1/r}(2l+1)^{(n-1)/r-1}$.

Given $\alpha$, choose an integer
$r$ so large that $\epsilon =1/(n-1)-\alpha-1/r$
is greater than $0$.
Then, with $m=(2l+1)^{n-1}$, we have
$m^\alpha (1-\delta (m,n))\leq (r!)^{1/r}/m^{\epsilon}$.
\end{proof}
\Rem{In the opposite direction, an example shows that,
for any real number $\alpha>2/(n-1)$,
$$
m^\alpha (1-\delta (m,n)) \to \infty
\text{\ as\ } m\to\infty\, .
$$
}
\begin{proof}
We shall use Lemma \ref{gen} with $r=0$.

Let $V=\{ (t_1,\ldots ,t_n)\in\Rr^n\st \sum_{i=1}^nt_i=1\}$
and let $\Delta\subseteq V$ be the simplex 
$\{ (t_i)\in V\st t_i\geq 0\}$.
The map $\pi : V\to S(\Rr^n)$, $v\mapsto v/\| v\|$ gives a
diffeomorphism from $V$ to an open subspace of the sphere
mapping geodesics (straight line segments) to geodesics (segments of 
great circles). Hence, there is a constant $c_n>0$ such that
$d(\pi (u),\pi (v))\geq c_n\| u-v\|$ for all $u,\, v$ in the compact
subpace $\Delta$.

Fix an integer $l\geq 1$ and let $W=\{ (t_i)\in\Delta \st lt_i\in\Zz
\text{\ for all $i$}\}$.
For $u,\, v\in W$, $u\not=v$, we see that $\| u-v\|^2\geq 2/l^2$.
List the vectors $\pi (v)$, $v\in W$, as $w_0,\ldots ,w_k$,
where $k+1=\binom{l+n-1}{n-1}\geq l^{n-1}/(n-1)!\,$.
So, for $i\not=j$, $d(w_i,w_j) \geq c_n\sqrt{2}/l$.

This shows that $\delta (k+1-n,n) \leq \cos (c_n\sqrt{2}/l)
=1 -c_n^2/l^2+\cdots$.
Finally, observe that $m=k+1-n\geq \frac{1}{2}l^{n-1}/(n-1)!$ if
$l$ is sufficiently large.
\end{proof}
\appendix
\Sect{On the cohomology of $\PO (\Rr^{2^r})$}
In this appendix we give an elementary proof of the result from
\cite[Lemma (8.1)]{proj} that we used in the proofs of Theorem 
\ref{simplex} and Lemma \ref{gen}. 
\Prop{\label{PSO}
Write $V=\Rr^{2^r n}$, where $r\geq 1$ and $n\geq 1$ is odd. 
Let $L$ be the
real line bundle over the projective orthogonal group
$\PO (V)$ associated with the double cover $\O (V)\to\PO (V)$.
Then $e(L)\in H^1(\PO (V);\,\Ff_2)$ satisfies
$e(L)^{2^r-1}\not =0$, $e(L)^{2^r}=0$.
Moreover, the restriction of $e(L)^{2^r-1}$ to the
projective special 
\hyphenation{ortho-gonal}
orthogonal group $\PSO (V)\leq \PO (V)$ is
non-zero.
}
\begin{proof}
Notice, first, that the vector bundle $L\otimes V$ is trivial:
it has a natural trivialization $L\otimes V \to V$ mapping
$g\otimes v$, where $g$ is an element of $\O (V)$, 
regarded as a generator of the
fibre of $L$ at $[g]\in\PO (V)$, and $v\in V$, to $gv$.
Hence $1= w(L\otimes V)=(1+e(L))^{2^rn} =(1+e(L)^{2^r})^n$ and,
since $n$ is odd, $e(L)^{2^r}=0$.

The diagonal map $\O (\Rr^{2^r}) \to (\O (\Rr^{2^r}))^n
\leq \O (\Rr^{2^rn})$ induces an inclusion of $\PO (\Rr^{2^r})$
in $\PO (V)$. To prove that $e(L)^{2^r-1}$ is non-zero it
is, therefore, enough to deal with the case $n=1$,
and for the remainder of the discussion we take
$V=\Rr^{2^r}$.

Given a Euclidean vector bundle $E$ over a base $B$,
we write $\O (\Rr^i,E)$ for the Stiefel bundle over
$B$ with fibre at $b\in B$ the Stiefel manifold
$\O (\Rr^i,E_b)$ of isometric linear maps
from $\Rr^i$ to the fibre $E_b$ of $E$ at $b\in B$.
 
Let $H$ be the Hopf bundle over the real projective
space $P(V)$ of $V$. Now there is a natural projection
$$
\pi :\O (V,H\otimes V) \to \PO (V)
$$
taking an isometric linear map $v\mapsto u\otimes g(v):
V\to \Rr u\otimes V$, where $g\in \O (V)$,
in the fibre over $[u]\in P(V)$ to $[g]\in\PO(V)$.
The pullback $\pi^*L$ is identified with the lift of $H$
(because $(-u)\otimes g(v)=u\otimes (-g)v$).

In terms of $x=e(H)\in H^1(P(V);\,\Ff_2)$, the total Stiefel-Whitney
class of $H\otimes V$ is equal to $(1+x)^{2^r}=1$, so that
$w_j(H\otimes V)=0$ for $j\geq 1$.

For each $i$, $0\leq i<2^r$, the fibration
$$
p_i :\O (\Rr^{i+1}, H\otimes V) \to \O (\Rr^i,H\otimes V)
$$
(restricting from $\Rr^{i+1}=\Rr^i\oplus\Rr$ to $\Rr^i$
and interpreted when $i=0$ as the projection
$p_0: S(H\otimes V) \to P(V)$)
is the sphere bundle $S(\zeta_i)$ of the complementary
$(2^{r}-i)$-dimensional vector bundle 
$\zeta_i$ over $\O (\Rr^i, H\otimes V)$.
Thus $\zeta_i\oplus\Rr^i$ is the pullback of $H\otimes V$,
and $e(\zeta_i)=w_{2^r-i}(\zeta_i)$ is the lift of 
$w_{2^{r}-i}(H\otimes V)$, which is zero.
Hence, from the Gysin sequence of the sphere bundle,
the induced homomorphism
$$
p_i^*:
H^*(\O (\Rr^i,H\otimes V);\, \Ff_2) \to
H^*(\O (\Rr^{i+1},H\otimes V);\, \Ff_2)
$$
is injective. 
In the top dimension $i=2^r-1$, when $\zeta_i$ is a line
bundle with $w_1=0$, $p_i$ is a trivial
bundle with fibre $S^0$ and $p_i^*$ is an injection
on each component. 

Hence, $(\pi^* e(L))^{2^r-1}$, which is the image
of $e(H)^{2^r-1}=x^{2^r-1}\in H^*(P(V);\,\Ff_2)$
under the composition 
$$
p_0\comp \cdots \comp p_{2^r-1} : \O (V,H\otimes V)
=\O (\Rr^{2^r},H\otimes V) \to
P(V)=\O (\Rr^0,H\otimes V),
$$
is non-zero on each of the two components.
\end{proof}
\Rem{Simpler proofs can be given in the
low dimensional cases $V=\Rr^{2^r}$, $r=1,\, 2,\, 3$.
In those dimensions there is a bilinear map
$\cdot : V\times V \to V$ satisfying $\| u\cdot v\|
=\| u\| \cdot \| v\|$ for $u,\, v\in V$,
given by complex, quaternionic or Cayley multiplication.
This determines a map $S(V) \to \O (V)$ taking $u\in S(V)$
to the orthogonal transformation $v\mapsto u\cdot v$, which 
induces an embedding of $P(V)$ into $\PO (V)$. 
The line bundle $L$ restricts to the Hopf bundle 
over $P(V)$, and so we see that $e(L)^{2^r-1}\not=0$.
}
\Rem{A result of Gitler and Handel \cite[Theorem 1.6]{gh} on projective
Stiefel manifolds can be established by the same method.
For $V=\Rr^{2^rn}$, with $n$ odd, consider, for $1\leq l\leq n$,
the projective Stiefel manifold
$\PO (\Rr^{2^rl},V)=\O (\Rr^{2^rl},V)/\{ \pm 1\}$ and the associated
real line bundle $L$. Then $e(L)^{2^rm-1}$ is non-zero if
$\binom{n}{j}$ is even for $n-m < j<l$,
and $e(L)^{2^rm}=0$ if $n-m<l$ and $\binom{n}{m}$ is odd.

For $e(\zeta_i)$ is the lift of $w_{2^rn-i}(H\otimes V)=
\binom{2^rn}{2^rn-i}x^{2^rn-i}$, which is $0$ if $i$ is not divisible
by $2^r$ and equal to $\binom{n}{n-j}x^{2^r(n-j)}$ if $i=2^rj$;
and the bundle $L\otimes V$ has a canonical trivial summand
$\Rr^{2^rl}$, so that $w_i(L\otimes V)=0$ for
$i>2^r(n-l)$.
}

\end{document}